# Bilateral Bounds for Norms of Solutions and Boundedness/Stability and Instability of Some Nonlinear Systems with Delays and Variable Coefficients

Mark A. Pinsky

**Abstract**. This paper presents a novel methodology for evaluating the boundedness, stability, and instability of some vector nonlinear systems with multiple time-varying delays and variable coefficients. The proposed technique develops two scalar counterparts for the initial vector system. The solutions to these scalar nonlinear equations, which also incorporate delays and variable coefficients, provide upper and lower bounds for the norms of solutions to the original vector equations with corresponding history functions. This enables evaluation of the dynamics of a vector system through the analysis of its scalar counterparts. This analysis can be accomplished using simplified analytical reasoning or straightforward simulations, which remain effective even for systems with a large number of coupled equations. Consequently, we introduced some novel boundedness, stability and instability criteria and estimated the radiuses of the balls containing the history functions which stem bounded/stabile solutions to the vector systems. Lastly, we validate our results in representative simulations that also assess their accuracy.

## 1. Introduction

Nonlinear vector systems with delays and variable coefficients play a pivotal role in analyzing the temporal dynamics of various models across biology, medicine, economics, physics, control science, and engineering. Comprehending the boundedness, stability, and instability of these systems is essential for grasping the complex temporal phenomena emerging in these diverse fields. The study of these systems not only deepens our understanding of these temporal behaviors but also enhances the predictive power and control of such models in real-world applications. Over the past several decades, a slew of monographs [5], [7], [8], [11]-[13], [15], [17], [19], [22], [24], [26] and research papers has emerged in this area, presenting a variety of techniques for assessing the stability of equilibrium solutions in nonlinear systems with delay. For additional insights and further references on this topic see review papers [2], [3], [14], [20] that provide overviews of the progress made in this field.

Significant advancement in the stability analysis of systems with delays has been marked by the pioneering work of N. Krasovskii [25] and B. Razumikhin [36], who extended the Lyapunov-function methodology to these systems. The merit of these approaches lies in the effectiveness of the Lyapunov-like functionals or functions they employ. However, finding appropriate functionals/functions for a broad class of vector nonlinear systems with delays remains a complex and challenging task. Lyapunov-like techniques applied to systems with delays often result in overly conservative stability criteria, a limitation that is particularly pronounced in systems with variable coefficients [28]. Nevertheless, for linear systems with delays and constant coefficients, these techniques lead to the development of computationally tractable method of linear matrix inequalities (LMI). For further details and comprehensive reviews on this subject, see [12], [14], [15], [23], [37].

The development of stability criteria for delay systems has also been influenced by Halanay's inequality [11], [12], [16], [29] and the Bohl-Perron theorem [5], [6], [13]. Despite these advancements, applying contemporary techniques to assess the stability of nonlinear systems with multiple delays remains challenging, often resulting in overly conservative criteria. Additionally, these methods rarely estimate the radii of the boundedness or stability balls that encompass the history functions stemming bounded or stable solutions, thereby reducing their practical utility. These difficulties are especially evident in vector nonlinear systems with variable delays [39].

A more challenging problem involves ensuring the boundedness of solutions to vector nonlinear systems with delays and external perturbations. This problem has been mainly explored within the framework of input-to-state stability (ISS), which is typically guaranteed under conservative assumptions regarding the stability of the corresponding homogeneous systems. For more details, see the recent survey [3] and additional references therein, as well as, e.g., [1] for a distinct approach to this problem.

The assessment of finite-time stability has also been extended to delay systems, as highlighted in a recent survey [38]. However, contemporary literature on the instability of delay systems remains somewhat limited. Some studies in this domain leverage Krasovskii's methodology [30], alongside other techniques [18], [27].

In [31]-[34] we bring forth a methodology that develops the scalar counterparts for a broad class of nonlinear systems of ODEs with variable coefficients. Solutions to these scalar differential equations provide upper bounds for the time histories of the norms of solutions to their multidimensional counterparts. This approach enables the assessment of the boundedness/stability characteristics of complex nonlinear systems by analyzing the dynamic behavior of their scalar complements. Additionally, it facilitates efficient estimations of the boundedness/stability regions for such multidimensional systems [34]. In [35] this approach was extended to certain vector nonlinear systems with delays, although this extension applies to a relatively narrow application domain.

This paper introduces a novel technique that develops a pair of scalar counterparts for a vector nonlinear system with multiple time-varying delays and variable coefficients. The solutions to these scalar nonlinear equations, which also incorporate delays and variable coefficients, provide upper and lower bounds for the norms of solutions to the corresponding vector systems with matched history functions. This allows for the assessment of the dynamics of such complex systems through simplified reasoning focused on their scalar counterparts, which can be efficiently analyzed via straightforward simulations or concise analytical reasoning. Consequently, we derive bilateral bounds for the time histories of the norms of solutions to a locally stable or locally bounded vector system with delays, as well as a lower bound for the norms of solutions to a locally unstable or unbounded system of this type.

As a result, we derive new criteria for the boundedness, stability, and instability of complex delay nonlinear systems and estimate the radii of the balls containing history functions that give rise to stable or bounded solutions for these systems. Finally, we validate our approach through representative simulations that also assess the accuracy of these techniques.

This current methodology, developed under an additional assumption applicable in various scenarios, has become more versatile and computationally tractable than our previous technique [35]. It appears to establish bilateral bounds for nonlinear systems involving multiple equations with variable delays and coefficients. However, in the common application domain, our earlier methodology generally provides more efficient estimates than the current one.

This paper is organized as follows. The next section presents the relevant notations, preliminary results, and the underlined equations. Section three outlines our methodology. Section four derives various criteria for the boundedness, stability, and instability of the systems under study based on the developed technique. Section five discusses the simulation results, and the last section concludes this study.

--------------------------------------------------------------------------------

Mark A. Pinsky, Depart. of Math. & Stat., University of Nevada. Reno, Reno NV 89557, USA, email: pinsky@unr.edu.

## 2. Notation, Mathematical Preliminaries, and Governing Equation

**2.1. Notation.** Firstly, we recall that symbols $\mathbb{R}$, $\mathbb{R}_{\geq 0}$, $\mathbb{R}_+$ and $\mathbb{R}^n$ stand for the sets of real, non-negative and positive real numbers, and real $n$-dimensional vectors, $\mathbb{N}$ is a set of real positive integers, $\mathbb{R}^{n \times n}$ is a set of $n \times n$-matrices, $I \in \mathbb{R}^{n \times n}$ is the identity matrix, $\mathbb{C}^n$ and $\mathbb{C}^{n \times n}$ are the sets of $n$-dimensional vectors and $n \times n$-matrices with complex entries. Next, $C([a,b]; \mathbb{R}^n)$, $C([a,b]; \mathbb{R}_+)$ and $C([a,b]; \mathbb{R}_{\geq 0})$ are the spaces of continuous real functions $\zeta : [a,b] \to \mathbb{R}^n$, $\zeta : [a,b] \to \mathbb{R}_+$ or $\zeta : [a,b] \to \mathbb{R}_{\geq 0}$, respectively, with the uniform norm $\|\zeta(t)\| := \sup_{t \in [a,b]} |\zeta(t)|$, where $|\cdot|$ stands for the Euclidean norm of a vector or the induced norm of a matrix, and $b$ can be infinity. Also note that $|x|_\infty = \sup_{i=1,\ldots,n}(|x_i|)$ and $|x|_1 = \sum_{i=1}^n |x_i|$, $\forall x \in \mathbb{R}^n$; and in the sequel of this paper the upper righthand derivative in $t$, $D^+ x(t) := \dot{x}(t)$.

**2.2. Preliminaries.** The comparison statements for solutions to scalar ODEs can be to scalar delay differential equations through the application of the method of steps, e.g., [35]. Next, we recall a statement from [35] that will be used in the subsequent sections of this paper.

**Lemma 1.** Let us write two scalar equations with delay in the following form,

$$\dot{u}_1 = f_1\big(t, u_1(t), u_1(t - h_1(t)), \ldots, u_1(t - h_m(t))\big), \ \forall t \in [t_0, t_*], \ u_1(t) = \varphi_1(t), \ \forall t \in [t_0 - \bar{h}, t_0]$$

$$\dot{u}_2 = f_2\big(t, u_2(t), u_2(t - h_1(t)), \ldots, u_2(t - h_m(t))\big), \ \forall t \in [t_0, t_*], \ u_2(t) = \varphi_2(t), \ \forall t \in [t_0 - \bar{h}, t_0]$$

where $\forall u_i \in \mathbb{R}$, continuous functions $f_i \in \mathbb{R}$ are locally Lipschitz in the second variables $\forall u_i \in J_+ \subset \mathbb{R}$, $\varphi_i(t) \in C\left(\left[t_0 - \bar{h}, t_0\right]; \mathbb{R}\right)$, $i = 1, 2$, and

$$h_i(t) \in C\left([t_0, \infty); \mathbb{R}_+\right), \max_i \sup_{\forall t \geq t_0} h_i(t) = \bar{h} < \infty, \min_i \inf_{\forall t \geq t_0} h_i(t) \geq \underline{h} > 0, i = 1, \ldots, m \quad (2.1)$$

Next, we assume that both of the above equations admit the unique solutions $u_i(t, \varphi) \in J_+$, $\forall t \in [t_0, t_*]$, where $t_*$ can be infinity, and that,

$$f_1(t, x_1, \ldots, x_{m+1}) \leq f_2(t, \bar{x}_1, \ldots, \bar{x}_{m+1}), \; x_1 \in \mathbb{R}, \; x_i \leq \bar{x}_i \in \mathbb{R}, \; i = 2, \ldots, m+1, \; \forall t \in [t_0, t_*],$$
$$\varphi_1(t) \leq \varphi_2(t), \; \forall t \in \left[t_0 - \bar{h}, t_0\right]$$

Then $u_1(t, \varphi_1) \leq u_2(t, \varphi_2)$, $\forall t \in [t_0, t_*]$.

In [35] we consider a vector nonlinear equation with delays and variable coefficients,

$$\dot{x} = B(t)x + f_*\left(t, x(t), x(t - h_1(t)), \ldots, x(t - h_m(t))\right) + F_*(t), \; \forall t \geq t_0$$
$$x(t) = \varphi(t), \; \forall t \in \left[t_0 - \bar{h}, t_0\right] \quad (2.2)$$

where $x \in \mathrm{N} \subset \mathbb{R}^n$, $0 \in \mathrm{N}$, a function $f_* \in \mathbb{R}^n$ is continuous in all variables and locally Lipschitz in the second one, $f_*(t, 0) = 0$, $\varphi(t) \in J \subset C\left(\left[t_0 - \bar{h}, t_0\right]; \mathbb{R}^n\right)$, $\|\varphi(t)\| := \sup|\varphi(t)|$, $\forall t \in \left[t_0 - \bar{h}, t_0\right]$, $\|\varphi(t)\| \leq \bar{\varphi} \in \mathbb{R}_{\geq 0}$, matrix $B(t) \in C\left([t_0, \infty); \mathbb{R}^{n \times n}\right)$, $F_*(t) = F_0 e(t)$, $e \in C\left([t_0, \infty); \mathbb{R}^n\right)$, $\|e(t)\| = 1$, $F_0 \in \mathbb{R}_{\geq 0}$, and functions $h_i(t)$ comply with (2.1). We also are going to use in the sequel the abridged notation, $x_*(t, \varphi) := x_*(t, t_0, \varphi)$, $\forall t \geq t_0$.

[35] derives a scalar complement to (2.2) that was written as follows,

$$\dot{u} = p(t)y + c(t)\left(L\left(t, u(t), u(t - h_1(t)), \ldots, u(t - h_m(t))\right) + |F(t)|\right)$$
$$u\left(t, |\varphi(t)|\right) = |\varphi(t)|, \; \forall t \in \left[t_0 - \bar{h}, t_0\right] \quad (2.3)$$

where $u(t) \in \mathbb{R}_{\geq 0}$, $p(t) = d\left(\ln|w(t)|\right)/dt$, $w(t)$ is the fundamental solution matrix for linear equation $\dot{x} = B(t)x$, $c(t) = |w(t)||w^{-1}(t)|$ is the running condition number of $w(t)$. Furthermore, a scalar function $L \in C\left([t_0, \infty) \times \mathbb{R}_{\geq 0}^{m+1}; \mathbb{R}_{\geq 0}\right)$, $L(t, 0) = 0$ sets up a nonlinear extension of the Lipschitz continuity condition as follows,

$$|f(t, \chi_1, \ldots, \chi_{m+1})| \leq L(t, |\chi_1|, \ldots, |\chi_{m+1}|), \; \chi_i \in \mathbb{R}^n, \; i = 2, \ldots, m+1,$$
$$\forall \chi = [\chi_1, \ldots, \chi_{m+1}]^T \in \Omega \in \mathbb{R}^{n(m+1)} \quad (2.4)$$

where $\Omega$ is a compact subset of $\mathbb{R}^{n(m+1)}$ containing zero and as prior, $f \in C\left([t_0, \infty) \times \mathbb{R}^{n(m+1)}; \mathbb{R}^n\right)$, see [35]. Note that Appendix1 also defines such scalar functions in a close form if $f(t, \chi_1, \ldots, \chi_{m+1})$ is a polynomial, power series or a rational function in $\chi_1, \ldots, \chi_{m+1}$, see also additional examples in Section 5 of this paper. Additionally, $L$ turns out to be a linear function in $|\chi_i|$ if $f$ is a linear function in the corresponding variables.

Note that in the remainder of this paper we assume for simplicity that $\Omega \equiv \mathbb{R}^{n(m+1)}$.

In turn, [35] infers that $|x(t, \varphi)| \leq u(t, |\varphi|)$, $\forall t \geq t_0$ which abridges analysis of boundedness/stability characteristics of (2.2).

Nonetheless, computing $|w(t)|$ can be awkward in high dimensions and $c(t)$ might approach infinity if $t \to \infty$ making the relevant inferences overconservative. In fact, let $B = diag\left(\lambda_1 \; \lambda_2\right)$, $\lambda_{1,2} \in \mathbb{R}$, $\lambda_1 > \lambda_2$. Then, $c(t) = e^{(\lambda_1 - \lambda_2)(t - t_0)}$ and $\lim_{t \to \infty} c(t) = \infty$.

To escape these limitations, we introduce a novel pair of scalar counterparts to equation (2.2) with $p = const$ and $c = 1$, incorporating bilateral bounds for $|x(t, \varphi)|$ that will facilitate our subsequent analysis.

## 2.3. Underlined Equation and Definitions.

Firstly, write (2.2) as follows,

$$\dot{x} = A_* x + G_*(t) x + f_*\left(t, x(t), x(t-h_1(t)),..., x(t-h_m(t))\right) + F_*(t), \ \forall t \geq t_0,$$
$$x(t) = \varphi(t), \ \forall t \in \left[t_0 - \bar{h}, t_0\right] \tag{2.5}$$

where $A_*(t_0) = \lim_{t \to \infty} \sup_{\forall t \geq t_0} (t-t_0)^{-1} \int_{t_0}^{t} B(s) ds \in \mathbb{R}^{n \times n}$ is a nonzero matrix of general position that assumes simple eigenvalues for all but possibly some isolated values of $t_0 \in \Upsilon \subset \mathbb{R}$ [4], and $G_*(t) = B(t) - A_*(t_0)$ is a zero mean matrix.

To simplify further referencing, we also acknowledge the homogeneous counterpart of (2.5),

$$\dot{x} = A_*(t) x + f_*\left(t, x(t), x(t-h_1(t)),..., x(t-h_m(t))\right), \ \forall t \geq t_0, \ x(t) = \varphi(t), \ \forall t \in \left[t_0 - \bar{h}, t_0\right] \tag{2.6}$$

and the following linear equation,

$$\dot{x} = A_*(t) x, \ \forall t \geq t_0, \ x(t_0, \varphi(t_0)) = \varphi(t_0) = x_0 \in \mathbb{R}^n \tag{2.7}$$

and assume that both (2.5) and (2.6) admit unique solutions $\forall \|\varphi(t)\| \leq \bar{\varphi}$, $\forall t \geq t_0 \in \Upsilon \subset \mathbb{R}$.

As is known, the solution to the last equation can be written as $x(t) = W(t, t_0) \varphi(t_0)$, where $W(t, t_0) = w(t) w^{-1}(t_0)$ is the transition (Cauchy) matrix and $w(t)$ is a fundamental solution matrix for (2.7).

Next, we lay down the standard definitions of the boundedness of solutions to (2.5) and stability of the trivial solution to (2.6) that are referenced in the sequel of this paper, see, e.g., [25] and [35].

**Definition 1**. Assume that $\varphi(t) \in J \subset C\left(\left[t_0 - \bar{h}, t_0\right]; \mathbb{R}^n\right)$ and (2.6) admits a unique solution $\forall \|\varphi(t)\| \leq \bar{\phi} > 0$. Then the trivial solution of equation (2.6) is called:

1) stable for the set value of $t_0$ if $\forall \varepsilon \in \mathbb{R}_+, \exists \delta_1(t_0, \varepsilon) \in \mathbb{R}_+$ such that $\forall \|\varphi\| < \delta_1(t_0, \varepsilon), |x(t, t_0, \varphi)| < \varepsilon, \forall t \geq t_0$. Otherwise, the trivial solution is called unstable.

2) uniformly stable if in the above definition $\delta_1(t_0, \varepsilon) = \delta_2(\varepsilon)$.

3) asymptotically stable if it is stable for given value of $t_0$ and $\exists \delta_3(t_0) \in \mathbb{R}_+$ such that $\lim_{t \to \infty} |x(t, t_0, \varphi)| = 0$, $\forall \|\varphi\| < \delta_3(t_0)$.

4) uniformly asymptotically stable if it is uniformly stable and in the previous definition $\delta_3(t_0) = \delta_4 = const$.

5) uniformly exponentially stable if $\exists \delta_5 \in \mathbb{R}_+$ and $\exists c_i \in \mathbb{R}_+, i=1,2$ such that,

$$|x(t, t_0, \varphi)| \leq c_1 \|\varphi\| \exp(-c_2(t-t_0)), \forall \|\varphi\| \leq \delta_5, \forall t \geq t_0$$

**Definition 2**. Let (2.5) admits a unique solution $\forall \|\varphi(t)\| \leq \bar{\phi}$, $\varphi(t) \in J \subset C\left(\left[t_0 - \bar{h}, t_0\right]; \mathbb{R}^n\right)$. A solution to equation (2.5) is called:

6) bounded for the set value $t_0$ if $\exists \delta_6(t_0) \in \mathbb{R}_+, \exists F_*(t_0) \in \mathbb{R}_{\geq 0}$ and $\exists \varepsilon_*(\delta_6, F_*) \in \mathbb{R}_+$ such that $|x(t, t_0, \varphi)| < \varepsilon_*$ $\forall t \geq t_0, \forall \|\varphi\| < \delta_6(t_0)$ and $\forall F_0 \leq F_*(t_0)$.

Otherwise, solution $x(t, t_0, \varphi)$ is called unbounded.

7) uniformly bounded if both $\delta_6(t_0) = \delta_7 \in \mathbb{R}_+$ and $F_*(t_0) = \hat{F}, \delta_7, \hat{F} \in \mathbb{R}_+$.

Next, we recall Krasovskii's definition of the stability under persistent perturbations [25, pp.161-164] that also can be called the robust stability. Consider a perturb equation to (2.6) in the following form,

$$\dot{x}_* = A_* x_* + G_*(t) x_* + f_*\left(t, x_*(t), x_*(t-h_1(t)),..., x_*(t-h_m(t))\right) +$$
$$R\left(t, x_*, x_*(t-h_1^*(t)),..., x_*(t-h_1^*(t))\right), \forall t \geq t_0 \tag{2.8}$$
$$x_*(t) = \varphi(t), \ \forall t \in \left[t_0 - \bar{h}, t_0\right]$$

where it is assumed that (I) function $f_*(t, \chi_1,...,\chi_{m+1})$ is Lipschitz continuous in $\chi_i \in \mathbb{R}^n$, $i=1,...,m+1$, (II) $f(t,0)=0$, (III) $R \in \mathbb{R}^n$ is a continuous function in all its variables but $R(t,0)$ might not be zero, and (IV) $h_i^* \in C([t_0,\infty);\mathbb{R}_+)$ comply with the conditions (2.4).

**Definition 3.** The trivial solution to (2.6) is called robustly stable if $\forall \varepsilon \in \mathbb{R}_+$, $\exists \Delta_i(\varepsilon) \in \mathbb{R}_+$, $i=1,2,3$ such that $|x_*(t,t_0,\varphi)|_\infty \leq \varepsilon$, $\forall t \geq t_0$ if $|R(t,\chi_1,...,\chi_{m+1})|_\infty < \Delta_1(\varepsilon)$, $\forall t \geq t_0$, $\forall |\chi_i|_\infty < \varepsilon$, $\|\varphi(t)\| < \Delta_2(\varepsilon)$, $|h_i(t) - h_i^*(t)| < \Delta_3(\varepsilon)$, $\forall t \geq t_0$, $i=1,...,m$, where $x(t,t_0,\varphi)$ is a solution to (2.8).

**Statement** (N. Krasovskii [25], p.162). Assume that conditions (I)-(IV) hold and that the trivial solution to (2.6) is uniformly asymptotically stable. Then it is robustly stable as well.

Lastly, we extend a definition of the finite-time stability (FTS) developed for ODEs [10] on the systems with delay.

**Definition 4.** Equations (2.5)/(2.6) are called

  a) finite-time stable (FTS) with respect to positive numbers $\eta_1, \eta_2, T \in \mathbb{R}_+$, $\eta_1 < \eta_2$ if the condition $\|\varphi\| < \eta_1$ implies that $|x(t,t_0,\varphi)| < \eta_2(\eta_1, T)$, $\forall t \in [t_0, t_0+T]$.

  b) finite-time contractively stable (FTCS) with respect to positive numbers $\eta_1, \eta_2, \eta_3, T \in \mathbb{R}_+$, $\eta_1 < \eta_3$ if it is FTS with respect of $\eta_1, \eta_2, T$ and $\exists t_1 \in (t_0, t_0+T)$ such that $|x(t,t_0,\varphi)| < \eta_3(\eta_1, T)$, $\forall t \in [t_1, t_0+T]$.

### 3. Delay Auxiliary Equations

This section derives a pair of scalar delay differential equations that provide upper and lower estimates for the norms of solutions to (2.5). We begin by outlining the comprehensive steps involved and then proceed to detail the underlying conditions.

Firstly, let us set that complex conjugate eigenvalues of $A_*(t_0)$, $\lambda_k(t_0) = \alpha_k(t_0) \pm i\beta_k(t_0)$, $i = \sqrt{-1}$, $1 \leq k \leq n_1 \leq n$, $n_1 \in \mathbb{N}$, and real eigenvalues of $A_*(t_0)$, $\lambda_k(t_0) = \alpha_k(t_0)$, $n_1 < k \leq n$, where $\alpha_k \in \mathbb{R}$ and $\beta_k \in \mathbb{R}_{>0}$. Additionally, we presume that $\alpha_k \geq \alpha_{k+1}$, $k \in [1, n-1]$ and also define a square diagonal matrices,

$$\alpha = diag(\alpha_1, \alpha_1, ..., \alpha_{n_1}, \alpha_{n_1}, \alpha_{n_1+1}, ..., \alpha_n), \beta = diag(\beta_1, -\beta_1, ..., \beta_{n_1}, -\beta_{n_1}, 0, ..., 0), \alpha, \beta \in \mathbb{R}^{n \times n}$$

**3.1** To derive the auxiliary equations with $p = const$ and $c = 1$, we rewrite (2.5) into the eigenbasis of $A_*$ as follows,

$$\dot{y} = Ay + G(t)y + f(t, y(t), y(t-h_1(t)),..., y(t-h_m(t))) + F(t),$$
$$y(t) = V^{-1}\varphi(t), \forall t \in [t_0 - \bar{h}, t_0] \tag{3.1}$$

where $y = V^{-1}x$, $y \in \mathbb{C}^n$, $V \in \mathbb{C}^{n \times n}$ is the eigenmatrix of $A_*$, $A = V^{-1}A_*V = \alpha + i\beta$, $G = V^{-1}G_*V$, $F(t) = V^{-1}F_*$, and $f(t, y(t), y(t-h_1(t)),..., y(t-h_m(t))) = V^{-1}f_*(t, Vy(t), Vy(t-h_1(t)),..., Vy(t-h_m(t)))$. As prior, we shall write the solution to the last equation as follows, $y(t,t_0,\varphi(t)) = y(t,\varphi)$.

Subsequently, we rewrite (3.1) as,

$$\dot{y} = (\lambda I + i\beta)y + (\alpha - \lambda I)y + G(t)y + f(t, y(t), y(t-h_1(t)),..., y(t-h_m(t))) + F(t)$$
$$y(t) = V^{-1}\varphi(t), \forall t \in [t_0 - \bar{h}, t_0] \tag{3.2}$$

where $\lambda \in \mathbb{R}$ is defined below. Then, we mention that the fundamental matrix of solutions for equation $\dot{y} = (\lambda I + i\beta)y$ can be written as follows,

$$w(t) = \exp(\lambda I + i\beta)(t - t_0) = e^{\lambda(t-t_0)}\exp(i\beta)(t-t_0) \tag{3.3}$$

Hence, $w(t_0) = w^{-1}(t_0) = I$ and $|w(t)| = e^{\lambda(t-t_0)}$, $|w^{-1}(t)| = e^{-\lambda(t-t_0)}$ since $|\exp i\beta(t-t_0)| = |\exp[-i\beta(t-t_0)]| = 1$.

Next, we write (3.2) in the integral form using the variation of parameters defined by matrix (3.3) as follows,

$$y(t,\varphi) = w(t)w^{-1}(t_0)V^{-1}\varphi(t_0) + w(t)\int_{t_0}^{t} w^{-1}(\tau)Q(\tau, y(\tau))d\tau, \quad \forall t \geq t_0,$$

$$y(t,\varphi) = V^{-1}\varphi(t), \quad \forall t \in [t_0 - \bar{h}, t_0]$$

where $Q(t, y(t), y(t-h_1(t)),..., y(t-h_m(t))) = (\alpha - \lambda I)y + G(t)y + f(t, y(t), y(t-h_1(t)),..., y(\tau-h_m(t))) + F(t)$.

Note that in subsequent formulas we set that $w^{-1}(t_0) = I$ due to (3.3).

Then, we admit the following relation,

$$|y(t,\varphi)| = \left| w(t)V^{-1}\varphi(t_0) + w(t)\int_{t_0}^{t} w^{-1}(\tau)Q(\tau, y(\tau), y(\tau-h_1(\tau)),..., y(\tau-h_m(\tau)))d\tau \right|, \quad \forall t \geq t_0, \quad (3.4)$$

$$|y(t,\varphi)| = |V^{-1}\varphi(t)|, \quad \forall t \in [t_0 - \bar{h}, t_0]$$

In turn, the applications of both the triangular inequality and its reverse counterpart yields the following pair of equations,

$$|y_1^{\pm}(t,\varphi)| = |w_{\pm}(t)V^{-1}\varphi(t_0)| \pm$$

$$|w_{\pm}(t)| \int_{t_0}^{t} |w_{\pm}^{-1}(\tau)| \left| Q(\tau, y_1^{\pm}(\tau), y_1^{\pm}(\tau-h_1(\tau)),..., y_1^{\pm}(\tau-h_m(\tau))) \right| d\tau, \quad \forall t \geq t_0, \quad (3.5)$$

$$|y_1^{\pm}(t,\varphi)| = |V^{-1}\varphi(t)|, \quad \forall t \in [t_0 - \bar{h}, t_0]$$

where $y_1^{\pm}(t,\varphi) \in \mathbb{C}^n$ and $w_{\pm}(t) = e^{\lambda_{\pm}(t-t_0)}\exp(i\beta)(t-t_0)$, $\lambda_{\pm} \in \mathbb{R}$, and the selections of subscripts/superscripts, i.e., $^+/_-$, are aligned with the selection of the alike signs in the right side of (3.5). Thus, (3.5) infers the following inequality $|y_1^{-}(t,\varphi)| \leq |y(t,\varphi)| \leq |y_1^{+}(t,\varphi)|, \quad \forall t \geq t_0$.

Next, to bring (3.5) into a tractable form, we apply to its right side the standard norm's inequalities and bound the nonlinear component via utility of (2.4). This prompts the following equations,

$$y_2^{\pm}(t,\varphi) = |w_{\pm}(t)V^{-1}\varphi(t_0)| \pm$$

$$|w_{\pm}(t)| \int_{t_0}^{t} |w_{\pm}^{-1}(\tau)| \left( q_{\pm}(\tau, y_2^{\pm}(\tau), y_2^{\pm}(\tau-h_1(\tau)),..., y_2^{\pm}(\tau-h_m(\tau))) \right) d\tau, \quad \forall t \geq t_0, \quad (3.6)$$

$$y_2^{\pm}(t,\varphi) = |V^{-1}\varphi|, \quad \forall t \in [t_0 - \bar{h}, t_0]$$

where

$$q_{\pm}(t, y_2^{\pm}(t), y_2^{\pm}(t-h_1(t)),..., y_2^{\pm}(t-h_m(t))) =$$

$$\pm\left( (|\alpha - \lambda_{\pm}I| + |G(t)|)y_2^{\pm} + L(t, y_2^{\pm}(t), y_2^{\pm}(t-h_1(t)),..., y_2^{\pm}(t-h_m(t))) + |F(t)| \right)$$

and $y_2^{\pm}(t) \in \mathbb{R}_{\geq 0}$. Note that $y_2^{+}(t) \geq 0, \forall t \geq t_0$ and $y_2^{-}(t)$ should be set at zero if $y_2^{-}(t) < 0$. Consequently, we conclude that $y_2^{-}(t,\varphi) \leq |y(t,\varphi)| \leq y_2^{+}(t,\varphi), \forall t \geq t_0$ since $|y_1^{+}(t,\varphi)| \leq y_2^{+}(t,\varphi), \forall \lambda_{+} \in \mathbb{R}$ and $|y_1^{-}(t,\varphi)| \geq y_2^{-}(t,\varphi), \forall \lambda_{-} \in \mathbb{R}$.

Next, we match the solutions to (3.6) with the solutions to the initial problem of the following pair of scalar differential equations,

$$\dot{y}_3^{\pm} = p_{\pm}(t)y_3^{\pm} \pm c(t)q_{\pm}(t, y_3^{\pm}), \quad \forall t \geq t_0$$

$$y_3^{\pm}(t, \phi) = \phi(t), \quad \forall t \in [t_0 - \bar{h}, t_0] \quad (3.7)$$

where $y_3^{\pm}(t) \in \mathbb{R}_{\geq 0}$, $p:[t_0, \infty) \to \mathbb{R}$, $c:[t_0, \infty) \to [1, \infty]$ and $\phi(t) \in C([t_0 - \bar{h}, t_0]; \mathbb{R}_{\geq 0})$.

Subsequently, we rewrite (3.7) in the integral form using the variation of parameters,

$$y_3^\pm(t) = e^{d_\pm(t,t_0)}\phi(t_0) \pm \int_{t_0}^t e^{-d_\pm(t,\tau)} c(\tau) q_\pm\left(\tau, y_3^\pm(\tau), y_3^\pm(\tau-h_1(\tau)),\ldots, y_3^\pm(\tau-h_m(\tau))\right) d\tau, \ \forall t \geq t_0 \quad (3.8)$$

$$y_3^\pm(t) = \phi(t), \ \forall t \in [t_0 - \overline{h}, t_0]$$

where $d_\pm(t,t_0) = \int_{t_0}^t p_\pm(s) ds$.

To determine $p_\pm(t)$ and $c(t)$, we match the right side of (3.6) and (3.8). Matching the first additions in these formulas, i.e., $|w_\pm(t) V^{-1} \varphi(t_0)|$ and $e^{d(t,t_0)} \phi(t_0)$, returns that,

$$|w_\pm(t) V^{-1} \varphi(t_0)| = \phi(t_0) \exp\left(\int_{t_0}^t p(s) ds\right)$$

Then, utility of (3.3) yields that $|w_\pm(t) V^{-1} \varphi(t_0)| = e^{\lambda_\pm(t-t_0)} |V^{-1}\varphi(t_0)|$ which implies that,

$$p_\pm = d \ln\left(|w_\pm(t) V^{-1}\varphi(t_0)|\right)/dt = \lambda_\pm, \ \forall \lambda_\pm \in \mathbb{R}$$

Note that the last relation yields that $d_\pm(t,\tau) = \lambda_\pm(t-\tau)$. To conclude that $c(t) \equiv 1$, we match the second editions on the right side of (3.6) and (3.8). This yields that

$$\int_{t_0}^t e^{-d(t,\tau)} c(\tau) q_\pm\left(\tau, y_3^\pm(\tau), y_3^\pm(\tau-h_1(\tau)),\ldots, y_3^\pm(\tau-h_m(\tau))\right) d\tau =$$

$$e^{\lambda t} \int_{t_0}^t e^{-\lambda\tau} c(\tau) q_\pm\left(\tau, y_3^\pm(\tau), y_3^\pm(\tau-h_1(\tau)),\ldots, y_3^\pm(\tau-h_m(\tau))\right) d\tau. \text{ Then, utility of (3.3) in (3.6) yields that } c(t) \equiv 1.$$

In turn, matching the history functions in (3.6) and (3.8), we resolve that $\phi(t) = |V^{-1}\varphi(t)|$.

Next, we write (3.7) as follows,

$$\dot{y}_3^\pm = (\lambda_\pm \pm |\alpha - \lambda_\pm I|) y_3^\pm \pm |G(t)| y_3^\pm \pm \left(L(t, y_3^\pm(t), y_3^\pm(t-h_1(t)),\ldots, y_3^\pm(t-h_m(t))) + |F(t)|\right) \quad (3.9)$$

$$y_3^\pm(t,\phi) = |V^{-1}\varphi(t)|, \ \forall t \in [t_0 - \overline{h}, t_0], \ \forall \lambda_\pm \in \mathbb{R}$$

and presume that both equations (3.9) admit unique solutions if $\|\varphi(t)\| \leq \overline{\varphi}$. This brings forth the following inequality,

$$y_3^-(t,\varphi) \leq |y(t,\varphi)| \leq y_3^+(t,\varphi), \ \forall t \geq t_0, \ \forall \lambda_\pm \in \mathbb{R} \quad (3.10)$$

**3.2** To develop a less conservative version of (3.9), we set that matrix, $D(t) := \text{Im}(diagG(t))$, $D \in \mathbb{R}^{n\times n}$, $\beta_*(t) = \beta + iD(t)$ and $g = G - iD$, and rewrite (3.2) as follows,

$$\dot{y} = (\lambda I + i\beta_*) y + (\alpha - \lambda I) y + g(t) y + f(t, y(t), y(t-h_1(t)),\ldots, y(t-h_m(t))) + F(t)$$

$$y(t) = V^{-1}\varphi(t), \ \forall t \in [t_0 - \overline{h}, t_0]$$

Then, the fundamental matrix of solution for equation $\dot{y} = (\lambda I + i\beta_*(t)) y$ takes the following form $w_*(t) = \exp(\lambda I + i\beta_*(t))(t-t_0)$, where $\beta_*$ is a diagonal matrix. This implies that $|w_*(t)| = e^{\lambda(t-t_0)}$, $|w_*^{-1}(t)| = e^{-\lambda(t-t_0)}$ which again attests that $c = 1$ and $p_\pm = \lambda_\pm$ and let us to cast (3.9) in a more efficient form.

**3.3** To further sharpen the estimates provided by (3.10) we select that $\lambda_+ := \min_{\lambda_+}(\lambda_+ + |\alpha - \lambda_+ I|)$ and $\lambda_- := \max_{\lambda_-}(\lambda_- - |\alpha - \lambda_- I|)$. Since matrix $\alpha - \lambda I$ is diagonal, $|\alpha - \lambda I| = \max_i |\alpha_i - \lambda|$, $i = 1,\ldots,n$ which yields that $\lambda_+ := \min_{\lambda_+}(\lambda_+ + \max_i |\alpha_i - \lambda_+|)$ and $\lambda_- := \max_{\lambda_-}(\lambda_- - \max_i |\alpha_i - \lambda_-|)$. Due to standard reasoning, these relations yield that $\lambda_+ = \alpha_1$ and $\lambda_- = \alpha_n$.

**3.4** Next, we enter (3.9) as the following pair of scalar equations,

$$\dot{Z} = (\alpha_1 + |g(t)|)Z + (L(t, Z(t), Z(t-h_1(t)), ..., Z(t-h_m(t))) + |F(t)|)$$
$$Z(t,\phi) = |V^{-1}\varphi(t)|, \forall t \in [t_0 - \bar{h}, t_0]$$
(3.11)

$$\dot{z} = (\alpha_n - |g(t)|)z - (L(t, z(t), z(t-h_1(t)), ..., z(t-h_m(t))) + |F(t)|)$$
$$z(t,\phi) = |V^{-1}\varphi(t)|, \forall t \in [t_0 - \bar{h}, t_0]$$
(3.12)

where $Z(t,\phi)$ and $z(t,\phi) \in \mathbb{R}_{\geq 0}$ since we shall set that $z(t,\phi) = 0$ if it takes negative values.

In turn, using that $y = V^{-1}x$, we streamline (3.10) as follows,

$$\frac{1}{|V^{-1}|} z(t, \phi(t)) \leq |x(t, \varphi(t))| \leq |V| Z(t, \phi(t)), \forall t \geq t_0$$
$$\phi(t) = |V^{-1}\varphi(t)|, \forall t \in [t_0 - \bar{h}, t_0]$$
(3.13)

Lastly, we reclaim the prior steps in the following statement.

**Theorem 1**. Assume that $f(t, \chi_1, ..., \chi_{m+1}) \in \mathbb{R}^n$, $\forall \chi_i \in \mathbb{R}^n$, $\forall t \geq t_0$ is a continuous function in all variables and locally Lipschitz in the second one, $f(t,0) = 0$, $\varphi(t) \in J \subset C([t_0 - \bar{h}, t_0]; \mathbb{R}^n)$, continuous scalar functions $h_i(t)$ are defined by (2.1), matrix $B(t) \in C([t_0, \infty); \mathbb{R}^{n \times n})$, $A_*(t_0) \in \mathbb{R}^{n \times n}$ is a matrix of general position that assumes simple eigenvalues for all but possibly some isolated values of $t_0 \in \Upsilon$, $F(t) = F_0 e(t)$, $e \in C([t_0, \infty); \mathbb{R}^n)$, $\|e(t)\| = 1$, $F_0 \in \mathbb{R}_{\geq 0}$, a scalar function $L(t, \xi_1, ..., \xi_{m+1}) \times \mathbb{R}_{\geq 0}^{m+1} \in C([t_0, \infty); \mathbb{R}_{\geq 0})$, $\xi_i \in \mathbb{R}_{\geq 0}$, $i = 1, ..., m+1$ is locally Lipschitz in the second variable, $L(t,0) = 0$, inequality (2.4) holds with $\Omega \equiv \mathbb{R}^{n(m+1)}$. Subsequently, we assume that equations (2.5), (2.6), (3.11) and (3.12) assume unique solutions $\forall \|\varphi(t)\| \leq \bar{\varphi}$, $\forall t \geq t_0$.

Then, inequality (3.13) holds for all but possibly some isolated values of $t_0 \in \Upsilon$, where $x(t, \varphi)$, $Z(t, |\varphi|)$ and $z(t, |\varphi|)$ are solutions to equations (3.1), (3.11) and (3.12) respectively.

**Proof**. In fact, under the above conditions the sequence of steps undertaken in subsections 3.1-3.4 of this section infers (3.11) □

4. **Boundedness, stability and instability of vector delay nonlinear time-varying equations**

To simplify further references, we lay down the homogeneous counterparts to (3.11) and (3.12) as follows,

$$\dot{Z} = (\alpha_1 + |g(t)|)Z + L(t, Z(t), Z(t-h_1(t)), ..., Z(t-h_m(t)))$$
$$Z(t,\phi) = |V^{-1}\varphi(t)|, \forall t \in [t_0 - \bar{h}, t_0]$$
(4.1)

$$\dot{z} = (\alpha_n - |g(t)|)z - L(t, z(t), z(t-h_1(t)), ..., z(t-h_m(t)))$$
$$z(t,\phi) = |V^{-1}\varphi(t)|, \forall t \in [t_0 - \bar{h}, t_0]$$
(4.2)

and subsequently assume that the last two equations admit unique solutions under the conditions of Theorem 1.

Next, we pretend that $r_i$ are the superior values of $\delta_i$ ($r_i = \sup \delta_i$) for which in Definition 1 the $i$-th statement, holds for solutions to (4.1) or in Definition 2, the statements 6 or 7 hold for solutions to (3.11).

Subsequently, we define the akin values $R_i = \sup \delta_i$ through the application of the prior ansatz to the initial equations (2.6) and (2.5), respectively, and, in turn, recall that $R_1 = R_1(t_0), R_3 = R_3(t_0), r_1 = r_1(t_0), r_3 = r_3(t_0)$, $R_6 = R_6(t_0, F_0)$, $r_6 = r_6(t_0, F_0), \forall F_0 \leq F_*(t_0)$, and $R_7 = R_7(F_0), r_7 = r_7(F_0), \forall F_0 \leq \hat{F}$.

Next, let $B_{R_i} := \varphi(t) \in C\left(\left[t_0 - \bar{h}, t_0\right], \mathbb{R}^n\right): \|\varphi\| \leq R_i$, $i = 1,...,7$ and $B_{r_i} := \varphi(t) \in C\left(\left[t_0 - \bar{h}, t_0\right], \mathbb{R}^n\right): \|V^{-1}\varphi\| \leq r_i$, $i = 1,...,7$ be the balls with radiuses $R_i$ and $r_i$, respectively, that are centered at the origin. Thus, $B_{R_i}$ encompasses the history functions stemming bounded or stable solutions for (2.5)/(2.6) whereas $B_{r_i}$ contains the history functions stemming akin solutions for (3.11)/(4.1). This leads to the following statements.

**Theorem 2**. Assume that the conditions of Theorem 1 are met and the trivial solution to equation (4.1) is either stable, uniformly stable, asymptotically stable, uniformly asymptotically stable, or exponentially stable. Then the trivial solution to equation (2.6) also is stable, uniformly stable, asymptotically stable, uniformly asymptotically stable, or exponentially stable, respectively.

Furthermore, under the conditions of 1- 5 of Definition 1, $B_{r_i} \subseteq B_{R_i}$, $i = 1,...,5$.

**Proof**. Both statements directly follow from the application of the rightmost part of inequality (3.13) to the solutions to equations (2.6) and (4.1) □

**Theorem 3**. Assume that the conditions of Theorem 1 are met and the trivial solution to equation (4.2) is unstable. Then the trivial solution to equation (2.6) is unstable as well. Furthermore, under the above conditions, a solution to vector equation (2.6) is bounded in norm from the below by the solution of (4.2) with matched history function.

**Proof**. Both statements directly follow from the application of the leftmost part of inequality (3.13) to the solutions to equations (2.6) and (4.1) □

**Theorem 4**. Assume that the conditions of Theorem 1 are met and that:

1. the solutions to the equation (3.11) are bounded $\forall \varphi(t): \|V^{-1}\varphi(t)\| \leq r_6$ for the set values of $t_0$, and $F_0 \leq F_*(t_0)$; or uniformly bounded for $F_0 \leq \hat{F}$ and $\forall \varphi(t): \|V^{-1}\varphi(t)\| \leq r_7$.

2. a solution to (3.12), $z(t, \phi)$, $\phi = |V^{-1}\varphi(t)|$ is unbounded.

Then,
1. The solutions to equation (2.5) with matched history functions are also either bounded or uniformly bounded for the same values of $t_0$ and $F_0$, respectively. Furthermore, under the prior conditions $B_{r_i} \subseteq B_{R_i}$, $i = 6,7$.

2. The solution to (2.5) with matched history function, i.e. $x(t, \varphi)$ is unbounded.

**Proof**. The proof of this statement directly follows from the application of inequality (3.13) to the solutions to equation (2.5) and its scalar counterparts (3.11)/(3.12) □

Thus, the radiuses of the boundedness/stability balls for the vector equations, (2.5)/ (2.6) can be estimated in simulations of the scalar equations (3.11)/(4.1). This task can be further abridged through utility of the following,

**Corollary 1**. Assume the conditions of Theorem 1 hold and let $Z_i(t, \phi_i(t))$ and $z_i(t, \phi_i(t))$

$\phi_i \in C\left(\left[t_0 - \bar{h}, t_0\right]; \mathbb{R}_{\geq 0}\right), \|\phi_i\| \leq \bar{\phi}$, $i = 1, 2$ are the solutions to (3.11) and (3.12), respectively, with $\phi_1 \leq \phi_2$, $\forall t \in \left[t_0 - \bar{h}, t_0\right]$. Then $Z_1(t, \phi_1) \leq Z_2(t, \phi_2)$ and $z_1(t, \phi_1) \leq z_2(t, \phi_2)$, $\forall t \geq t_0$.

**Proof**. In fact, the solutions to equations (3.11) and (3.12) cannot intersect in $Z \times t$ or $z \times t$ planes due to uniqueness of solutions to these equations presumed in Theorem 1 □

Note that the above statement also follows from Lemma 1 under the assumption that $L(t, \xi_1,...,\xi_{m+1})$ is a nondecreasing function in all its variables starting from the second one that can be naturally assured if $f(t, \chi_1,..., \chi_{m+1})$, $\chi_i \in \mathbb{R}^n$ is a polynomial/power series in $(\chi_1,..., \chi_{m+1})$ with time-dependent coefficients, see Appendix 1.

**Corollary 2**. Assume that the conditions of Theorem 1 hold, $Z_1(t,c)$, $z_1(t,c)$, $c \in \mathbb{R}_{\geq 0}$, $\forall t \in [t_0 - \bar{h}, t_0]$ and $Z_2(t,\phi)$, $z_2(t,\phi)$, $\phi \in C([t_0 - \bar{h}, t_0]; \mathbb{R}_{\geq 0})$, $\|\phi(t)\| \leq c \leq \bar{\phi}$, $\forall t \in [t_0 - \bar{h}, t_0]$ are the solutions to (3.11)/(3.12), respectively, with constant and variable history functions. Then, $Z_2(t,\phi) \leq Z_1(t,c)$ and $z_2(t,\phi) \leq z_1(t,c)$, $\forall t \geq t_0$

**Proof**. The proof of this statement immediately follows from Corollary 1 □

Thus, the values of $r_i$ can be seamlessly approximated in simulations of the corresponding scalar delay equations over a set of constant and positive history functions.

In turn, we show that the robust stability condition for the trivial solution of a vector equation (2.6) can be attested through the alike assumptions made on its scalar counterpart.

Firstly, we lay down a scalar auxiliary equation to (2.6) as follows,

$$\dot{Z}_* = (\alpha_1 + |g(t)|)Z_* + L(t, Z_*(t), Z_*(t - h_1(t)),...,Z_*(t - h_m(t))) + L_R(t, Z_*(t), Z_*(t - h_1^*(t)),...,Z_*(t - h_m^*(t)))$$
$$Z_*(t,\phi) = |V^{-1}\varphi(t)|, \forall t \in [t_0 - \bar{h}, t_0] \quad (4.3)$$

where $Z_* \in \mathbb{R}_{\geq 0}$, $L$ is defined prior, $L_R \in \mathbb{R}_{\geq 0}$ is a continuous function in all its variables, $L_R(t,0)$ might not be equals zero, and

$$|R(t, \chi_1,...,\chi_{m+1})| \leq L_R(t, |\chi_1|,...,|\chi_{m+1}|), \chi_i \in \mathbb{R}^n, i = 2,...,m+1, \forall \chi = [\chi_1,...,\chi_{m+1}]^T \in \Omega \in \mathbb{R}^{n(m+1)} \quad (4.4)$$

This embraces the following,

**Theorem 5**. Assume that the conditions of Theorem 1 hold, $L_R \in \mathbb{R}_{\geq 0}$ is a continuous function in all its variables, inequality (4.4) holds, conditions (I)-(IV) of the Statement hold for equation (4.3), and the trivial solution to (4.1) is uniformly asymptotically stable. Then the trivial solution to (2.6) is robustly stable.

**Proof**. In fact, under the above conditions, the trivial solution to (4.3) is robustly stable, i.e. $Z_*(t, |\varphi|) < \varepsilon$, $\forall t \geq t_0$ if $|R(t, \chi_1,...,\chi_{m+1})|_\infty < \Delta_1(\varepsilon)$, $\forall t \geq t_0$, $\forall |\chi_i|_\infty < \varepsilon$, $\|\varphi(t)\| < \Delta_2(\varepsilon)$, $|h_i(t) - h_i^*(t)| < \Delta_3(\varepsilon)$, $\forall t \geq t_0$, $i = 1,...,m$, see Definition 3. Then, due to (4.4), (4.3) is the scalar counterpart to equation (2.8). Thus, the application of (3.13) to solutions of the last two equations yields that $|x_*(t, \varphi(t))| \leq |V| Z_*(t, \phi(t))$ and, in turn, assures this statement □

Next, we show that the application of FTS/FTCS conditions to scalar equation (3.11)/(3.10) ensures the relevant stability conditions for their vector counterparts (2.5)/(2.6).

**Theorem 6**. Assume that:
  a) Equations (3.11)/(4.1) are FTS with respect $\eta_1, \eta_2, T \in \mathbb{R}_+$, $\eta_1 < \eta_2$. Then equations (2.5)/(2.6) are FTS with respect $\eta_1, \eta_2, T$ as well.
  b) Equation (3.11)/(4.1) are FTCS with respect to $\eta_1, \eta_2, \eta_3, T \in \mathbb{R}_+$, $\eta_1 < \eta_3$. Then equations (2.5)/(2.6) are FTCS with respect $\eta_1, \eta_2, \eta_3, T$ as well.

**Proof**. In fact, both above statements are assured by (3.11) □

Note that equations (3.11)/(312) can be further reduced to their autonomous counterpart if $\sup_{t \geq t_0} L(t, \xi_1,...,\xi_{m+1}) = L_s(\xi_1,...,\xi_{m+1}) < \infty$, $\forall \xi_i \in \mathbb{R}_{\geq 0}$, $\sup_{t \geq t_0} g(t) = g_s < \infty$, and $\sup_{t \geq t_0} F(t) = F_s < \infty$. Clearly, the latter conditions are frequently met in applications and the solutions to these autonomous equations bound from the above and below the solutions to (3.11)/(3.12), respectively, and, thus, can be assimilated in (3.13) as well.

Furthermore, applications of the Lipschitz-like condition,

$$L(t, \xi_1,...,\xi_{m+1}) \leq \sum_{i=1}^{m+1} \mu_i^+(t, \bar{\xi})\xi_i, \forall t \geq t_0, \xi_i \in \mathbb{R}_{\geq 0}, \xi = [\xi_1,...,\xi_{m+1}]^T, \forall |\xi| \leq \bar{\xi} > 0$$, to the right sides of (3.11) and (3.12) yield the scalar linear equations with solutions bounding from the above and below the matching solutions to (3.11)/(3.12). The latter equations can be turned into autonomous linear delay scalar equations under the conditions of the prior ansatz, which abets their subsequent analysis but mounts the conservatism of subsequent estimates.

## 5. Simulations

In this section, the developed methodology is applied to simulate the bilateral bounds for the time histories of the norms of solutions to two systems, each coupling two Van der Pol-like (VDP) or Duffing-like (DUF) oscillators with delays, time-varying nonperiodic coefficients, external perturbations, and either typical dissipative or conservative nonlinear components. These types of systems, when delays are omitted, are commonly found in various applications [21]. The simulations presented in this section utilize MATLAB's DDE23 code and sample results from a much broader set of simulations that validate the key inequality (3.11) and support our subsequent conclusions.

### 5.1. Simulated Equations

Equation (2.5) can be turned into such systems if we assume that,

$$A_* = \begin{pmatrix} 0 & 1 & 0 & 0 \\ -(\omega_1^2 + d) & -\chi_1 & d & 0 \\ 0 & 0 & 0 & 1 \\ d & 0 & -(\omega_2^2 + d) & -\chi_2 \end{pmatrix}, G_* = \begin{pmatrix} 0 & 0 & 0 & 0 \\ -g_{21}(t) & 0 & 0 & 0 \\ 0 & 0 & 0 & 0 \\ 0 & 0 & -g_{43}(t) & 0 \end{pmatrix},$$

$$F_*(t) = \begin{pmatrix} 0 & F_{*1}(t) & 0 & F_{*2}(t) \end{pmatrix}^T$$

where $\omega_1^2 = 1$, $\omega_2^2 = 4$, $d = 4$, $F_{*i} = F_{0i} \sin q_i t$, $i = 1, 2$, $q_1 = 5.43$, $q_2 = 10$, $g_{21} = a_1 \sin r_1 t + a_2 \sin r_2 t$, $g_{43} = b_1 \sin s_1 t + b_2 \sin s_2 t$, $r_1 = 3.14$, $r_2 = 6.15$, $s_1 = 3.1$, $s_2 = 6.28$. The values of other parameters used in specific simulations are given in Appendix 2.

Next, we assume that for VDP system,

$$f_* = G_*(t) x(t - h_1) + \begin{pmatrix} 0 & -\mu_1 x_2^3(t - h_2) & 0 & -\mu_2 x_4^3(t - h_2) \end{pmatrix}^T$$

and for DUF system,

$$f_* = G_*(t) x(t - h_1) + \begin{pmatrix} 0 & -\mu_1 x_1^3(t - h_2) & 0 & -\mu_2 x_3^3(t - h_2) \end{pmatrix}^T$$

where $x \in \mathbb{R}^4$, and $h_1$ and $h_2 \in \mathbb{R}_+$. Next, let $V = [v_{ij}]$, $i, j = 1, ..., 4$ is the eigenvector matrix for matrix $A_*$, $A = diag(\lambda_1, \lambda_1^*, \lambda_2, \lambda_2^*)$, $\lambda_k = \alpha_k + i\beta_k$, $\lambda_k^* = \alpha_k - i\beta_k$, $k = 1, 2$, $\alpha_1 \geq \alpha_2$, as prior $G(t) = V^{-1} G_* V$, $F(t) = V^{-1} F_*(t)$, and for VDP system,

$$f = G(t) y(t - h_1) + V^{-1} \begin{pmatrix} 0 & -\mu_1 \left( \sum_{k=1}^{4} v_{2k} y_k(t - h_2) \right)^3 & 0 & -\mu_2 \left( \sum_{k=1}^{4} v_{4k} y_k(t - h_2) \right)^3 \end{pmatrix}^T \quad (5.1)$$

, and for DUF system, $f = G(t) y(t - h_1) + V^{-1} \begin{pmatrix} 0 & -\mu_1 \left( \sum_{k=1}^{4} v_{1k} y_k(t - h_2) \right)^3 & 0 & -\mu_2 \left( \sum_{k=1}^{4} v_{3k} y_k(t - h_2) \right)^3 \end{pmatrix}^T$.

Subsequently, for both of our sample systems we derive that,

$$L(t, z) = |G(t)| z(t - h_1) + \|V^{-1}\| abs(\mu_1 \kappa_1 + \mu_2 \kappa_2) z^3(t - h_2) \quad (5.2)$$

, where for VDP system, $\kappa_1 = \left( \sum_{k=1}^{4} abs(v_{2k}) \right)^3$, $\kappa_2 = \left( \sum_{k=1}^{4} abs(v_{4k}) \right)^3$, for DUF system $\kappa_1 = \left( \sum_{k=1}^{4} abs(v_{1k}) \right)^3$, $\kappa_2 = \left( \sum_{k=1}^{4} abs(v_{3k}) \right)^3$, and $abs(q) = |q|$ if $Im(q) = 0$ and $abs(q) = \sqrt{a^2 + b^2}$ if $q = a + ib$.

Let us recap the derivation of (5.2) for VDP system since its derivation for DUF system follows the same steps. Let $Y_i = \sum_{k=1}^{4} v_{ik} y_k$, $i = 2, 4$ and $Y = \begin{pmatrix} 0 & -\mu_1 Y_2^3 & 0 & -\mu_2 Y_4^3 \end{pmatrix}^T$. Then (5.1) takes the form, $f = G(t) y + V^{-1} Y$.

Subsequently, the application of standard inequalities, $|V^{-1} Y| \leq |V^{-1}||Y|$, $|Y|_2 \leq |Y|_1 = abs(\mu_1 Y_2^3) + abs(\mu_1 Y_4^3)$ and $|y_k^m| \leq |y|^m$, $m \in \mathbb{N}$, returns that $|Y|_1 = abs\left( \mu_1 \left( \sum_{k=1}^{4} v_{2k} y_k \right)^3 \right) + abs\left( \mu_2 \left( \sum_{k=1}^{4} v_{4k} y_k \right)^3 \right) = (abs(\mu_1) \kappa_1 + abs(\mu_2) \kappa_2) |y|^3$

which yields lastly that $|V^{-1} Y| \leq |V^{-1}| (abs(\mu_1) \kappa_1 + abs(\mu_2) \kappa_2) z^3(t - h_2)$.

Thus, equations (4.1) and (4.2) for both VDP and DUF systems can be written as follows,

$$\dot{Z} = (\alpha_1 + |g(t)|)Z + |G_*(t)|Z(t-h_1) +$$
$$|V^{-1}|\left(\left(abs(\mu_1)\kappa_1^3 + abs(\mu_2)\kappa_2^3\right)Z^3(t-h_2) + |F(t)|\right) \quad (5.3)$$
$$Z(t,\phi) = |V^{-1}\varphi(t)|, \forall t \in [t_0 - \bar{h}, t_0]$$
$$\dot{z} = (\alpha_2 - |g(t)|)z - |G_*(t)|z(t-h_1) -$$
$$|V^{-1}|\left(\left(abs(\mu_1)\kappa_1^3 + abs(\mu_2)\kappa_2^3\right)z^3(t-h_2) + |F(t)|\right) \quad (5.4)$$
$$z(t,\phi) = |V^{-1}\varphi(t)|, \forall t \in [t_0 - \bar{h}, t_0]$$

where $\kappa_1$ and $\kappa_2$ take different values for VDP and DUF systems as was indicated prior.

Note that simulations of both VDP and DUF systems adopt the same history function, $\varphi(t) = [x_{0,1},...,x_{0,4}]^T \cos t$, $x_{0,1} \in \mathbb{R}_+$ and $x_{0,i} = 0, i \neq 1$.

### 5.2. Simulation of VDP and DUF systems

The results of simulations for locally stable VDP and DUF systems are shown in Figure 1, where solid, dashed, and dash-dotted lines represent the time histories of the norms of solutions to the initial equations (2.5) or (2.6), and their scalar counterparts (5.3) or (5.4), respectively, simulated with matched history functions. All our simulations yield qualitatively similar results for both of our sample systems. Note that the sets of parameters for featured simulations are resumed in Appendix 2.

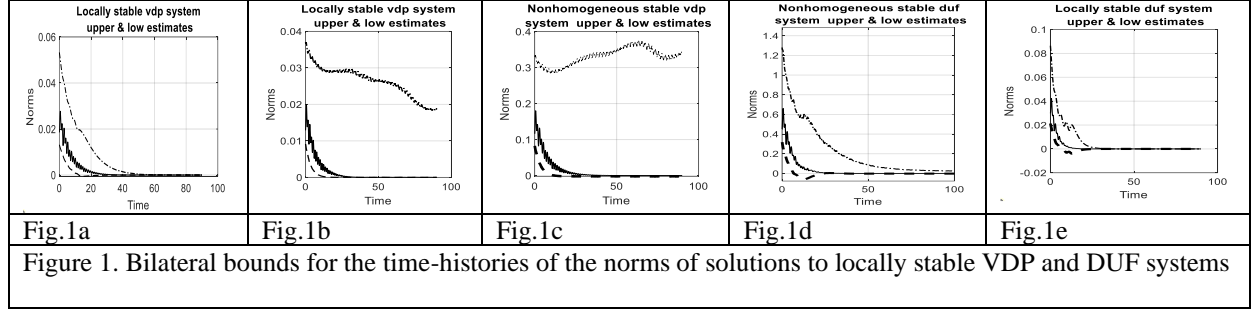

| Fig.1a | Fig.1b | Fig.1c | Fig.1d | Fig.1e |

Figure 1. Bilateral bounds for the time-histories of the norms of solutions to locally stable VDP and DUF systems

We observe that upper estimates frequently prove to be more conservative compared to their lower counterparts when calculated for the same set of parameters. This tendency is also evident when comparing the norms of solutions for the nonhomogeneous equation (2.5) to those of the homogeneous equation (2.6). Specifically, the upper estimates for the nonhomogeneous case are consistently more conservative than those for the homogeneous case, as illustrated in Figure 1c.

In both scenarios, the conservatism of the upper bound increases as the parameter $x_{0,1}$ grows closer to the threshold value at which the solution to equation (5.3) becomes unstable or unbounded. Yet, even when the solution to equation (5.3) becomes unstable the corresponding solution to equation (2.6) might remain stable, or even asymptotically stable.

The stability and boundedness properties of both the initial vector equations and their scalar counterparts are influenced by the magnitudes of the time-varying and friction components. However, the scalar equations are typically more sensitive to variation of these components, as is expected. For example, in Figure 1a, $a_i = b_i = 0$, while in Figure 1b, $a_i = b_i = 0.1$. As a result, the estimates derived from Figure 1a are more efficient than those from Figure 1b.

Furthermore, modifying the durations of delays appears to have a relatively minor impact on the accuracy of these estimates.

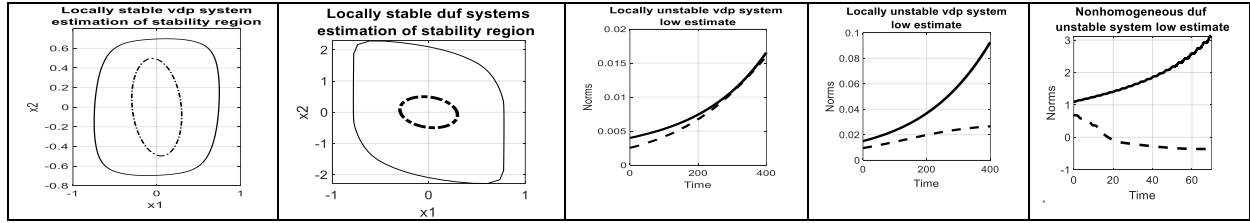

| Fig 2. Estimation of stability regions for locally stable VDP and DUF systems | Figure 3. Lower bounds for the time-histories of the norms of solutions to locally unstable VDP and DUF systems |
|---|---|

To further validate our approximations, we assume that $\varphi(t) \equiv x_0 = const$, $x_0 \in \mathbb{R}^4$ and estimate the boundaries of boundedness/stability regions for vector equations (2.5) or (2.6) and their scalar counterparts. These regions, i.e., $\vartheta_i \subset \mathbb{R}^4$, $i = 1, 2$ consist of zero vector and all other vectors $x_0$ that stem bounded/stable solutions to (2.5)/(2.6), $i = 1$, or akin solutions to their scalar counterparts, $i = 2$. Due to (3.11), $\vartheta_2 \subseteq \vartheta_1$ and the running time for estimating the boundary of $\vartheta_1$ scales as $m^n$, where $n$ is the number of independent variables in (2.5)/(2.6) and $m$ is the number of the discretization points for each of these variables.

Thus, to streamline such simulations, we estimate the projections of $\vartheta_1$ onto $x_1 \times x_2$ coordinate plane. For this purpose, we express $x_0$ in double polar coordinates, i.e., $r_i, \theta_i$, $r_i \in \mathbb{R}_{\geq 0}$, $\theta_i \in [0, 2\pi]$, $i = 1, 2$, discretize $\theta_i$ with step $\pi/60$ and for each pair of $\theta_i$ find the minimal values of $r_i$ for which $|x(t, x_0)|$ rapidly surpasses the threshold value that is set at $10^5$. A similar approach was applied to (5.3) to determine the minimal value of $|x_0| = |x_0|_{min}$ that implies alike behavior of $Z(t, |x_0|)$. Following this, we estimate the boundary of stability/boundedness region using the equation $\|V^{-1} x_0\| = |x_0|_{min}$, see the inner curves on the plots in Figure 2. Notably, the running time for such simulations is practically insensitive to changes in $n$ and $m$. This time can be further reduced since $Z(t, |x_0|)$ monotonically increase in $|x_0|$, $\forall t$.

Our approach provides seamless, though rather conservative, estimates of the boundaries of the stability regions for both the VDP and DUF systems. Nevertheless, these estimates validate our theoretical inferences, as shown in Figure 2, where the coordinates of the inner curves have been scaled by 10 for enhanced clarity.

Figure 3 presents samples from our simulations, illustrating the lower bounds for locally unstable VDP and DUF systems. As anticipated, our estimates tend to be more conservative for nonhomogeneous systems, as demonstrated in Figure 3c. Despite this conservatism, our approach often effectively identifies unstable behavior and provides reasonable estimates for the norms of solutions in locally unstable homogeneous systems, as is seen in Figures 3a and 3b.

In conclusion, all the simulations we conducted support our theoretical inferences.

## 6. Conclusion and Future Research

This paper develops a novel approach to analyzing the boundedness/stability and instability of some vector nonlinear systems with multiple time-varying delays and variable coefficients, which are common in various applications. To achieve this, we derive two scalar equations whose solutions provide upper and lower bounds for the time-histories of the norms of solutions to nonlinear vector systems that include delays and variable coefficients. This allows for the evaluation of the dynamics of vector systems by analyzing the solutions to their scalar counterparts, which can be further supported by seamless simulations. The running times of such simulations remain virtually unaffected by the system's complexity.

Furthermore, this approach aggregates the parameters and history functions, thereby simplifying the evaluation of sensitivity and robustness properties of complex delay systems.

As a result, we derived novel criteria for boundedness, stability, and instability, and estimated the radii of the balls enclosing the history-functions which stem bounded and stable solutions for a broad class of vector delay and nonlinear systems.

The outcomes of this effort were authenticated on representative simulations involving two systems with typical dissipative and conservative nonlinear components, variable coefficients, and distinct delays. The simulations confirm our inferences and measure their effectiveness.

Future research will aim to improve the accuracy of our estimates by combining the current technique with recursive approximations. Additionally, we intend to target functional differential equations with neutral and distributed delays and relevant control applications.

**Appendix 1**. Let us show how to derive inequality (2.4) using an example that can be naturally extended onto more complex cases. Assume that $x = [x_1 \ x_2]^T$ and a vector-function $f$ is defined, e.g., as follows, $f = \begin{bmatrix} a_1(t)x_1 x_2^3(t-h_1) & a_2(t)x_1^2(t-h_2) \end{bmatrix}^T$. Then $|f|_2 \leq |f|_1 \leq |a_1||x_1||x_2^3(t-h_1)| + |a_2||x_1^2(t-h_2)|$
$\leq |a_1||x(t)||x(t-h_1)|^3 + |a_2||x(t-h_2)|^2$, where we use that $|x_i^n| \leq |x|^n$, $|x_i^n(t-h)| \leq |x(t-h)|^n$, $i=1,2, n \in \mathbb{N}$. Clearly, this inference can be extended to the power series and some rational functions. In fact, assume again that $x \in \mathbb{R}^2$ and $g(x) = 1/(1+x_2(t-h))$, then $|g| \geq 1/(1-|x(t-h)|)$. The utility of this bound normally would require that $1-|x(t-h)| \neq 0, \forall t \geq t_0$.

**Appendix 2.** Figure 1a, $\mu_i = -3$, $a_i = b_i = 0$, $\chi_i = 0.4$, $i=1,2$, $h_1 = 10, h_2 = 12$. Figure 1b, $\mu_i = -0.5$, $a_i = b_i = 0.1$, $i=1,2$, $\chi_1 = 0.2$, $\chi_2 = 0.4$, $h_1 = 1, h_2 = 2$. Figure 1c, $\mu_i = -0.01$, $a_i = 0$, $b_i = 0.1$, $i=1,2$, $\chi_1 = 0.2$, $\chi_2 = 0.4$, $h_1 = 10, h_2 = 12, F_{0,1} = 0, F_{0,2} = 0.001$, Figure 1d, $\mu_i = -0.01$, $a_i = 0$, $b_i = 0.1$, $\chi_i = 0.6$, $i=1,2$, $h_1 = 10, h_2 = 12, F_{0,1} = 0, F_{0,2} = 0.001$. Figure 1e, $\mu_i = -0.01$, $a_i = 0$, $b_i = 0.1$, $\chi_i = 0.6 i=1,2, h_1 = 10, h_2 = 12$, $h_1 = 10, h_2 = 12$, Figure 2a, $\mu_i = -3$, $a_i = b_i = 0$, $\chi_i = 0.6$, $i=1,2$, $h_1 = 10, h_2 = 12$. Figure 2c, $\mu_i = -3$, $a_i = b_i = 0$, $\chi_i = 0.6$, $i=1,2$, $h_1 = 10, h_2 = 12$. Figure 3a, $\mu_i = -3$, $a_i = b_i = 0$, $\chi_i = 0.4$, $i=1,2$, $h_1 = 10, h_2 = 12$. Figure 3b, $\mu_i = -0.3$, $a_i = b_i = 0$, $i=1,2$, $\chi_1 = 0.2$, $\chi_2 = 0.4$, $h_1 = 10, h_2 = 12$. Figure 3c, $\mu_i = -0.01$, $a_i = b_i = 0$, $i=1,2$, $\chi_1 = 0.2$, $\chi_2 = 0.4$, $h_1 = 10, h_2 = 12, F_{0,1} = 0, F_{0,2} = 0.001$.

**Acknowledgement**. The simulation code for this study was developed by Steve Koblik.
Data sharing is not applicable to this article as no datasets were generated or analyzed during the current study.
The author declares that he has no conflict of interest.